\def\titlep{Generalized permutative representations of Cuntz algebras}
\documentclass[11pt]{article}
\usepackage{graphicx,ifthen}

\font\germ=eufm10 at12pt

\def\goth#1{\hbox{\germ#1}}

\newcommand{\qed}{\hbox{\rule[-2pt]{3pt}{6pt}}}
\newcommand{\qedh}{\hfill\qed \\}
\newcommand{\vv}{\vspace{.3in}}
\newcommand{\vep}{\varepsilon}

\newtheorem{Thm}{Theorem}[section]

\newtheorem{ex}[Thm]{Example}
\newtheorem{defi}[Thm]{Definition}
\newtheorem{lem}[Thm]{Lemma}
\newtheorem{conj}[Thm]{Conjecture}
\newtheorem{prop}[Thm]{Proposition}
\newtheorem{cor}[Thm]{Corollary}

\newcommand{\ww}{\vv\noindent}

\def\sc#1{S({\bf C}^{#1})}

\def\sct#1#2{S({\bf C}^{#1})^{\otimes #2}}

\def\sci{S({\bf C}^{N})^{\infty}}

\def\cal#1{\mathcal #1}
\def\con{{\cal O}_{N}}
\def\edot{=1,\ldots,N}
\def\pr{{\it Proof.}\quad}
\def\zen{z=(z^{(n)})}
\def\evp{eventually periodic}
\def\asp{asymptotically periodic}

\def\co#1{{\cal O}_{#1}}
\def\ltn{l_{2}({\bf N})}

\def\mattwo#1#2#3#4{
\left(
\begin{array}{cc}
#1&#2\\
#3&#4\\
\end{array}
\right)}
\addtocounter{footnote}{1}
\def\sftt#1{
\setcounter{equation}{0}
\addtocounter{footnote}{1}
\section{#1}
}

\def\ssft#1{\subsection{#1}}

\def\cls{\quad\clearpage}

\begin{document}
%
% Personal data
%
\def\autherp{Katsunori Kawamura}
\def\emailp{e-mail:kawamura@kurims.kyoto-u.ac.jp.}
\title{{\bf \titlep}}
\author{\autherp\footnote{\emailp}\\
{\it \addressp}}
%\date{}
\def\inputcls#1{\cls\input #1.txt}
\def\plan#1#2{\par\noindent\makebox[.5in][c]{#1}
\makebox[.1in][l]{$|$}
\makebox[3in][l]{#2}\\}
\def\nset#1{\{1,\ldots,N\}^{#1}}
\def\bfsnl{{\rm BFS}_{N}(\Lambda)}
\def\blank{\cls \quad}
\newcommand{\mline}{\noindent
\thicklines
\setlength{\unitlength}{.08mm}
\begin{picture}(1000,5)
\put(0,0){\line(1,0){1600}}
\end{picture}
\par
}
\def\enda{{\rm End}{\cal A}}
\def\endcon{{\rm End}\con}
\def\scn{S({\bf C}^{N})^{\otimes n}}
\newcommand{\ctk}{S({\bf C}^{N})^{\otimes k}}
\newcommand{\ct}{S({\bf C}^{N})}
%\begin{frontmatter}

\def\scp{S({\bf C}^{N})^{\otimes p}}
\def\scm#1{S({\bf C}^{N})^{\otimes #1}}

%
%%%%%%%%%%%%%%%%%%%%%%%%%%%%%%%%%%%%%%%%%%
%\input comm.txt
%%%%%%%%%%%%%%%%%%%%%%%%%%%%%%%%%%%%%%%%%%
\pagestyle{plain}
\setcounter{page}{1}
\setcounter{section}{0}
\setcounter{footnote}{0}

%
% Personal data
\def\addressp{College of Science and Engineering Ritsumeikan University,
1-1-1 Noji Higashi, Kusatsu, Shiga 525-8577,Japan
}

\title{{\bf \titlep}}
\author{\autherp\footnote{\emailp}\\
{\it \addressp}}
\begin{center}
{\Large \titlep}

\ww
\autherp\footnote{\emailp}\\
{\it \addressp}\\
\end{center}

%\maketitle

%%%%%%%%%%%%%%%%%%%%%%%%%%%%%%%%%%%
%
% Abstract
%
\begin{abstract}
We introduce representations of the Cuntz algebra $\con$ which are 
parameterized by sequences in the set of unit vectors in ${\bf C}^{N}$. 
These representations are natural generalizations of 
permutative representations by Bratteli-Jorgensen and Davidson-Pitts.
We show their existence, irreducibility, equivalence,
uniqueness of irreducible decomposition and decomposition formulae
by using parameters of representations.
\end{abstract}
%
%\begin{center}
%\begin{keyword}
%Cuntz algebra\sep representation
%{\it MSC:} 47D25
%\end{keyword}
%\end{center}
%\end{frontmatter}

%%%%%%%%%%%%%%%%%%%%%%%%%%%%%%%%%%%%%%%
%
%  Section 1
%
\sftt{Introduction}
\label{section:first}
Let $\con$ be the Cuntz algebra with canonical generators
$s_{1},\ldots,s_{N}$ and
let $({\cal H},\pi)$ be a representation
of $\con$ with a cyclic vector $\Omega$ which satisfies
%
% Equation 1.1
%
\begin{equation}
\label{eqn:introduction}
\pi(s_{1})\Omega=\Omega.
\end{equation}
Then $({\cal H},\pi)$ is irreducible and unique up to unitary equivalence,
and it is equivalent to a permutative representation by \cite{BJ,DaPi2,DaPi3}.
The restriction $({\cal H},\pi|_{\con^{U(1)}})$ 
on the fixed-point subalgebra $\con^{U(1)}$ 
with respect to the $U(1)$-gauge action is also irreducible.
Especially, $\co{2}^{U(1)}$ is isomorphic to the CAR algebra ${\cal A}$ 
with canonical generators $\{a_{n}:n\geq 1\}$ which satisfy
$a_{n}a_{m}^{*}+a_{m}^{*}a_{n}=\delta_{nm}I$ and
$a_{n}a_{m}+a_{m}a_{n}=0$ for each $n,m\geq 1$.
Define an embedding $\varphi$ of ${\cal A}$ into $\co{2}$ by
\[\varphi(a_{1})\equiv s_{1}s_{2}^{*},\quad
\varphi(a_{n})\equiv \sum_{J\in\{1,2\}^{n-1}}
s_{J}s_{1}s_{2}^{*}\beta_{2}(s_{J}^{*})\quad(n\geq 2)\]
where $s_{J}\equiv s_{j_{1}}\cdots s_{j_{n-1}}$
for $J=(j_{1},\ldots,j_{n-1})$ and
$\beta_{2}$ is an automorphism of $\co{2}$ defined by
$\beta_{2}(s_{1})\equiv s_{1}$ and $\beta_{2}(s_{2})\equiv -s_{2}$. 
Then $\varphi({\cal A})=\co{2}^{U(1)}$.
Further $({\cal H},\pi\circ \varphi)$ 
is equivalent to the Fock representation of ${\cal A}$ with the vacuum $\Omega$
with respect to the annihilator $\{a_{n}\}_{n\geq 1}$, that is,
\[({\cal H},\pi\circ \varphi)\sim Fock,\quad
\pi(\varphi(a_{n}))\Omega=0\quad(n\geq 1).\]
According to such facts, there are many applications for representations
of the CAR algebra from representations of $\con$ \cite{AK1,AK3,IWF01}.
Further permutative representations distinguish
unitary equivalence classes of a certain class of endomorphisms of $\con$ 
effectively by their branching laws \cite{PE01,PE02}.

By generalizing (\ref{eqn:introduction})
for a unit vector $z=(z_{1},\ldots,z_{N})\in {\bf C}^{N}$,
we define a representation 
$({\cal H},\pi)$ with a cyclic vector $\Omega$ such that
%
% Equation 1.2
% 
\begin{equation}
\label{eqn:gpequation}
\pi(z_{1}s_{1}+\cdots+z_{N}s_{N})\Omega=\Omega.
\end{equation}
We denote such representation by $GP(z)$.
Then
$GP(z)$ exists uniquely up to unitary equivalence and irreducible for each $z$.
Further $GP(z)$ and $GP(y)$ are equivalent if and only if $z=y$.
The equation (\ref{eqn:gpequation}) is 
equivalent to an eigenequation of the Perron-Frobenius operator 
for a representation of $\con$ associated with a dynamical system \cite{PFO}.
In this way, these representations have many applications and
their studies are mathematically interesting.

In this article, we generalize such representations as follows:
Define $\sc{N}\equiv \{z\in {\bf C}^{N}:\|z\|=1\}$,
$\sct{N}{k}\equiv \{z^{(1)}\otimes \cdots
\otimes z^{(k)}:z^{(i)}\in\sc{N},\,i=1,\ldots,k\}$
for $k\geq 1$ and $\sc{N}^{\infty}\equiv 
\{(z^{(n)})_{n\in {\bf N}}:z^{(n)}\in \sc{N}\}$.
For $z=(z_{i})_{i=1}^{N}\in {\bf C}^{N}$,
define $s(z)\equiv z_{1}s_{1}+\cdots+z_{N}s_{N}$.
%
% Definition 1.1
%
\begin{defi}
\label{defi:gpdef}
Let $({\cal H},\pi)$ be a representation of $\con$.
\begin{enumerate}
%(i)
\item
For $z=z^{(1)}\otimes \cdots
\otimes z^{(k)}\in \sct{N}{k}$, $({\cal H},\pi)$ is $GP(z)$ 
if there is a unit cyclic vector $\Omega\in {\cal H}$ such that 
$\pi(s(z))\Omega=\Omega$ and 
$\{\pi(s(z^{(i)})\cdots s(z^{(k)}))\Omega\}_{i=1}^{k}$
is an orthonormal family in ${\cal H}$
where $s(z)\equiv s(z^{(1)})\cdots s(z^{(k)})$.
%(ii)
\item
For $z=(z^{(n)})_{n\in {\bf N}}\in\sc{N}^{\infty}$,
$({\cal H},\pi)$ is $GP(z)$ 
if there is a unit cyclic vector $\Omega\in {\cal H}$ such that
$\{\pi(s(z^{(n)})^{*}\cdots s(z^{(1)})^{*})\Omega:n\in {\bf N}\}$
is an orthonormal family in ${\cal H}$.
\end{enumerate}
For both cases, we call $\Omega$ the {\it GP vector} of $({\cal H},\pi)$.
We call $({\cal H},\pi)$  in (i) ({\it resp}. (ii))
the GP (=generalized permutative) representation 
with a cycle ({\it resp}. a chain).
\end{defi}

\noindent
In $\S$ \ref{section:third}, we show that the case (ii) is equivalent to 
an induced product representation in \cite{BC}.

Define $\sct{N}{*}\equiv\coprod_{k\geq 1}\ctk$ and 
$\sc{N}^{\#}\equiv\sct{N}{*}\cup \sci$.
%
% Theorem 1.2
%
\begin{Thm}
\label{Thm:mainfirst}
For any $z\in \sc{N}^{\#}$, 
there is a representation $({\cal H},\pi)$ of $\con$ which is $GP(z)$.
Such representation is unique up to unitary equivalence.
\end{Thm}

\noindent
By Theorem \ref{Thm:mainfirst}, the symbol $GP(z)$ makes sense as both
a representation and an equivalence class of representations of $\con$.

We define three periodicities on $\sc{N}^{\#}$ as follows.
For $z\in \sct{N}{*}$, $z$ is {\it periodic} if there are $y\in \sct{N}{l}$
and $p\geq 2$ such that $z$ is equal to the tensor power $y^{\otimes p}$
of $y$.
An element $z=(z^{(n)})_{n\in {\bf N}}$ in $\sci$ is {\it \evp}
if there are positive integers $p$, $M$ and a sequence
$(c_{n})_{n\geq M}$ in $U(1)\equiv \{a\in {\bf C}:|a|=1\}$
such that $z^{(n+p)}=c_{n}z^{(n)}$ for each $n\geq M$.
$z=(z^{(n)})_{n\in {\bf N}}$ is {\it \asp}
if there is a positive integer $p$ such that 
$\sum_{n=1}^{\infty}(1-|<z^{(n)}|z^{(n+p)}>|)<\infty$.
The eventual periodicity is properly stronger than the asymptotic periodicity.
%
% Theorem 1.3
% Irreducibility
%
\begin{Thm}
\label{Thm:mainsecond}
\begin{enumerate}
%(i)
\item
For $z\in \sct{N}{*}$, 
$GP(z)$ is irreducible if and only if $z$ is nonperiodic.
%(ii)
\item
For $z\in \sci$, 
$GP(z)$ is irreducible if and only if $z$ is non\asp.
\end{enumerate}
\end{Thm}

We define an equivalence relation $\sim$ among $\sc{N}^{\#}$.
For $z=z^{(1)}\otimes \cdots \otimes z^{(k)}\in\sct{N}{k}$ and 
$y\in \sct{N}{l}$, $z\sim y$
if $l=k$ and $y=z^{(\sigma(1))}\otimes \cdots\otimes z^{(\sigma(k))}$
for some $\sigma\in {\bf Z}_{k}$.
For $z=(z^{(n)})_{n\in {\bf N}}$ and $y=(y^{(n)})_{n\in {\bf N}}$ in $\sci$,
$z\sim y$ if there are nonnegative integers $p,q$
such that 
$\sum_{n=1}^{\infty}(1-|<z^{(n+p)}|y^{(n+q)}>|)<\infty$.
For each $z\in \sct{N}{*}$ and $y\in\sci$, we define $z\not\sim y$.
For two representations $({\cal H},\pi)$ and $({\cal H}^{'},\pi^{'})$,
$({\cal H},\pi)\sim({\cal H}^{'},\pi^{'})$ means that $({\cal H},\pi)$
and $({\cal H}^{'},\pi^{'})$ are unitarily equivalent.

%
% Theorem 1.4
% Equivalence
%
\begin{Thm}
\label{Thm:mainthird}
For $z, y\in \sc{N}^{\#}$,
$GP(z)\sim GP(y)$ if and only if $z\sim y$.
\end{Thm}

\noindent
By Theorem \ref{Thm:mainfirst}, \ref{Thm:mainsecond} and \ref{Thm:mainthird},
the irreducible decomposition of GP representation makes sense.

%
% Theorem 1.5
% Decomposition
%
\begin{Thm}
\label{Thm:mainfourth}
\begin{enumerate}
%(i)
\item
If $p\geq 2$ and $z\in \sct{N}{*}$ is nonperiodic,
then the following irreducible decomposition holds:
\[GP(z^{\otimes p})
\sim  GP(\zeta_{1}\cdot z)\oplus\cdots\oplus GP(\zeta_{p}\cdot z).\]
where $\zeta_{j}\equiv e^{2\pi\sqrt{-1}(j-1)/p}$.
This decomposition is unique up to unitary equivalence and multiplicity free.
%(ii)
\item
If $z\in \sci$ is \evp,
then there is a nonperiodic element $y\in \sct{N}{*}$
such that the following direct integral decomposition holds:
%
% Equation 1.3
%
\begin{equation}
\label{eqn:gpchain}
GP(z)\sim \int_{U(1)}^{\oplus}GP(c\cdot y)\, d\eta(c)
\end{equation}
where $\eta$ is the Haar measure of $U(1)$.
If there is another $y^{'}\in \sct{N}{*}$ which 
satisfies (\ref{eqn:gpchain}) with respect to $z$,
then there is $c_{0}\in U(1)$ such that $y^{'}\sim c_{0}\cdot y$.
\end{enumerate}
\end{Thm}

\noindent
By Theorem \ref{Thm:mainfourth} (ii),
it is understood that the reason why cycles and chains are simultaneously
treated in this article.

In $\S$ \ref{section:second}, we show Theorem \ref{Thm:mainfirst}.
In $\S$ \ref{section:third},
we show states of the Cuntz algebras associated with GP representations
and relations between GP representations with a chain
and induced product representations.
In $\S$ \ref{subsection:fourthone}, 
we prove a part of Theorem \ref{Thm:mainsecond}.
In $\S$ \ref{subsection:fourthtwo}, we prove Theorem \ref{Thm:mainthird}.
In $\S$ \ref{subsection:fourththree},
we prove the remaining part of Theorem \ref{Thm:mainsecond}
and Theorem \ref{Thm:mainfourth}.
In $\S$ \ref{section:fifth},
we explain that any cyclic permutative representation is a GP representation
and we introduce examples which are neither permutative representation
nor the rotation of permutative representation by the canonical action of $U(N)$.
At last, we discuss a remaining problem in $\S$ \ref{subsection:fifththree}.

%%%%%%%%%%%%%%%%%%%%%%%%%%%%%%%%%%
%
% Section 2
%
\sftt{Existence, basis and uniqueness}
\label{section:second}
For $N\geq 2$, let $\con$ be the {\it Cuntz algebra} \cite{C}, that is,
a C$^{*}$-algebra which is universally generated by
generators $s_{1},\ldots,s_{N}$ satisfying
$s_{i}^{*}s_{j}=\delta_{ij}I$ for $i,j\edot$ and
$s_{1}s_{1}^{*}+\cdots +s_{N}s_{N}^{*}=I$.
Let $\alpha$ be the canonical $U(N)$ action on $\con$.
In this article, any representation means a unital $*$-representation.
%
% Proposition 2.1
%
\begin{prop}
\label{prop:existence}
For each $z\in \sc{N}^{\#}$, 
there is a representation of $\con$ which is $GP(z)$.
\end{prop}
%
% Proof
%
\pr
(i)
Let $z=z^{(1)}\otimes \cdots\otimes z^{(k)}\in \sct{N}{k}$ and
$z^{(l)}=(z_{1}^{(l)},\ldots,z^{(l)}_{N})$ for $l=1,\ldots,k$.
Define $Y_{k}\equiv \{1,\ldots,k\}\times {\bf N}$ and a representation
$\pi_{0}$ of $\con$ on ${\cal H}\equiv l_{2}(Y_{k})$ by
$\pi_{0}(s_{i})e_{l,m}\equiv e_{l-1,N(m-1)+i}$ for $i\edot$ and $(l,m)\in Y_{k}$
where ${\bf N}\equiv \{1,2,3,\ldots\}$,
$\{e_{x}\}_{x\in Y_{k}}$ is the canonical basis of  ${\cal H}$ and
we define $e_{0,m}\equiv e_{k,m}$ for $m\in {\bf N}$.
Because $e_{l,1}=\pi_{0}(s_{1}^{l-1})^{*}e_{1,1}$ for each $l=2,\ldots,k$,
$\{e_{x}\}_{x\in Y_{k}}\subset\pi_{0}(\con)e_{1,1}$.
This implies that $({\cal H},\pi_{0})$ is cyclic.
Choose a sequence $(g(l))_{l=1}^{k}$ in $U(N)$ which satisfies
$g(l)_{j1}=z^{(l)}_{j}$ for each $l=1,\ldots,k$ and $j\edot$
where $g(l)=(g(l)_{ij})_{i,j=1}^{N}$.
Define a new representation $\pi$ on ${\cal H}$ by
$\pi(s_{i})e_{l,m}\equiv\pi_{0}( \alpha_{g(l-1)^{*}}(s_{i}))e_{l,m}$
for $i\edot$ and $(l,m)\in Y_{k}$
where we define $g(0)\equiv g(k)$.
Since $\alpha_{g(l)}(s_{1})=s(z^{(l)})$,
we have $\pi(s(z^{(l-1)}))e_{l,1}=e_{l-1,1}$ for $l=1,\ldots,k$
where we define $z^{(0)}\equiv z^{(k)}$.
From this, $\pi(s(z))e_{1,1}=e_{1,1}$.
Because $\{e_{x}\}_{x\in Y_{k}}\subset\pi(\con)e_{1,1}$,
$({\cal H},\pi)$ is cyclic.
Therefore $({\cal H},\pi)$ is $GP(z)$.

\noindent
(ii)
Let $z=(z^{(n)})_{n\in {\bf N}}\in \sci$,
$z^{(n)}=(z^{(n)}_{1},\ldots,z^{(n)}_{N})\in\ct$
and $Y_{\infty}\equiv {\bf Z}\times {\bf N}$.
Define a representation $\pi_{0}$ of $\con$ on 
${\cal H}\equiv l_{2}(Y_{\infty})$ by
$\pi_{0}(s_{i})e_{n,m}\equiv e_{n-1,\,N(m-1)+i}$ for $i\edot$ 
and $(n,m)\in Y_{\infty}$.
Then $\{\pi_{0}((s_{1}^{*})^{n})e_{0,1}\}_{n\in {\bf N}}
=\{e_{n,1}\}_{n\in {\bf N}}$.
For convenience,
we extend $z$ as $z=(z^{(n)})_{n\in{\bf Z}}$ by
$z^{(-n)}\equiv (1,0,\ldots,0)\in\sc{N}$ for each $n\geq 0$.
Choose a sequence $(g(n))_{n\in {\bf Z}}$ in $U(N)$ which satisfies
%
% Equation 2.1
%
\begin{equation}
\label{eqn:secondd}
g(n)_{j1}=z^{(n-1)}_{j}\quad(j\edot,\,n\in {\bf Z}).
\end{equation}
Define a new representation $\pi$ of $\con$ on ${\cal H}$
by $\pi(s_{i})e_{n,m}\equiv\pi_{0}(\alpha_{g(n)^{*}}(s_{i}))e_{n,m}$
for $(n,m)\in Y_{\infty}$.
Because $s(z^{(n-1)})=\alpha_{g(n)}(s_{1})$,
$\{\pi(s(z[n]))^{*}e_{0,1}\}_{n\in{\bf N}}=\{e_{n,1}\}_{n\in {\bf N}}$
is an orthonormal family in ${\cal H}$
where $z[n]\equiv z^{(1)}\otimes \cdots \otimes z^{(n)}$.
Because $\{e_{x}\}_{x\in Y_{\infty}}\subset\pi(\con)e_{0,1}$,
$({\cal H},\pi)$ is cyclic.
Hence $({\cal H},\pi)$ is $GP(z)$.
\qedh

In order to show the uniqueness of a GP representation,
we construct a basis of the representation space of
a given GP representation which is 
unique up to unitary equivalence.

For $z=z^{(1)}\otimes\cdots\otimes z^{(k)}\in S({\bf C}^{N})^{\otimes k}$,
let $({\cal H},\pi)$ be $GP(z)$ with the GP vector $\Omega$.
Define a unit vector
%
% Equation 2.2
%
\begin{equation}
\label{eqn:bbs}
E_{n}\equiv\pi(s(z^{(n)})\cdots s(z^{(k)}))\Omega\quad(n=1,\ldots,k).
\end{equation}
Then $\pi(s(z^{(n-1)}))E_{n}=E_{n-1}$ for $n=2,\ldots,k$
and $\pi(s(z^{(k)}))E_{1}=E_{k}$.
We see that $<E_{n}|E_{n^{'}}>=\delta_{nn^{'}}$ for $n,n^{'}=1,\ldots,k$.
Choose an orthonormal family $\{z^{(n,j)}\}_{j=1}^{N}$ in ${\bf C}^{N}$
such that $z^{(n,1)}=z^{(n)}$, and define
\[E_{n,j}\equiv \pi(s(z^{(n,j)}))E_{n+1}\quad(j\edot,\,n=1,\ldots,k)\]
where $E_{k+1}\equiv E_{1}$.
Then $E_{n,1}=E_{n}$ and $\{E_{n,j}:(n,j)\in\{1,\ldots,k\}\times \nset{}\}$
is an orthonormal family in ${\cal H}$.
Let $\nset{0}\equiv \{0\}$ and $\nset{*}$ be the union of $\nset{m}$ 
with respect to every $m=0,1,2,\ldots$
Define
\[E_{J,n,j}\equiv \pi(s_{J})E_{n,j}\quad(J\in\nset{*},\,n=1,\ldots,k,\,j=2,\ldots,N)\]
where we define $s_{0}\equiv I$ for convenience.
Then 
\[{\cal B}(z)\equiv \{E_{n}\}_{n=1}^{k}\cup
\{E_{J,n,j}:(J,n,j)\in\nset{*}\times \{1,\ldots,k\}\times \{2,\ldots,N\}\}\]
is also an orthonormal family in ${\cal H}$.
Define a subset $\Lambda(z)$ of $\sct{N}{*}$ by
\[\mbox{$\Lambda(z)\equiv\coprod_{m\geq 0}
\{\Lambda^{(m)}_{1}(z)\sqcup\cdots\sqcup \Lambda^{(m)}_{k}(z)\}$}\]
where $\Lambda^{(0)}_{n}(z)\equiv 
\{z^{(n)}\otimes \cdots \otimes z^{(k)}\}$ for $n=1,\ldots,k$,
$\Lambda^{(1)}_{1}(z)\equiv$ $\{z^{(k,j)}\}_{j=2}^{N}$,
$\Lambda^{(1)}_{n}(z)\equiv 
\{z^{(n-1,j)}\otimes z^{(n)}\otimes\cdots \otimes z^{(k)}\}_{j=2}^{N}$
for $n=2,\ldots,k$,
$\Lambda^{(m)}_{n}(z)\equiv 
\{\vep_{J}\otimes x:x\in \Lambda^{(1)}_{n}(z),\,J\in \{1,\ldots,N\}^{m-1}\}$
for $n=1,\ldots,k$ and $m\geq 2$ where
$\{\vep_{i}\}_{i=1}^{N}$ is the canonical basis of ${\bf C}^{N}$
and $\vep_{J}\equiv\vep_{j_{1}}\otimes\cdots \otimes\vep_{j_{m}}$
when $J=(j_{1},\ldots,j_{m})\in\{1,\ldots,N\}^{m}$.

%
% Proposition 2.2
%
\begin{prop}
\label{prop:main}
For $z\in \sc{N}^{\#}$,
let $({\cal H},\pi)$ be $GP(z)$ with the GP vector $\Omega$.
Define a family ${\cal C}_{z}\equiv \{E_{x}\in {\cal H}:x\in \Lambda(z)\}$
of unit vectors in ${\cal H}$ by
%
% Equation 2.3
%
\begin{equation}
\label{eqn:basis}
E_{x}\equiv \pi(s(x))\Omega\quad\quad(x\in \Lambda(z)).
\end{equation}
Then ${\cal C}_{z}$ is a complete orthonormal basis of ${\cal H}$.
\end{prop}
%
% Proof
%
\pr
We can verify that ${\cal C}_{z}={\cal B}(z)$.
Hence ${\cal C}_{z}$ is an orthonormal family in ${\cal H}$.
On the other hand, $V_{z}\equiv{\rm Lin}\langle\{E_{n}\}_{n=1}^{k}\rangle$
satisfies $\pi(s_{J}^{*})V_{z}\subset V_{z}$ for each $J\in \nset{*}$.
Furthermore we see that 
$\overline{\pi(\con) V_{z}}={\cal H}_{0}\equiv\overline{
{\rm Lin}\langle\{\pi(s_{J})\Omega:J\in \nset{*}\}\rangle}$
where the overline means the closure in ${\cal H}$.
Because $\Omega$ is a cyclic vector, ${\cal H}={\cal H}_{0}$.
By definition of $\Lambda(z)$, 
${\rm Lin}\langle{\cal C}_{z}\rangle$ is dense in ${\cal H}_{0}$.
In consequence, the statement holds.
\qedh

%%%%%%%%%%%%%%%%%%%%%%%%%%%%%%%%%%%

Fix $\zen_{n\in {\bf N}}\in\sci$.
We extend as $z=(z^{(n)})_{n\in{\bf Z}}$ by $z^{(-n)}\equiv \vep_{1}$
for $n\geq 0$.
Choose $(g(n))_{n\in{\bf Z}}$ which satisfies (\ref{eqn:secondd})
and $g(-n)=I$ for $n\geq 0$.
Let $({\cal H},\pi)$ be $GP(z)$ with the GP vector $\Omega$.
Define
%
% Equation 2.4
%
\begin{equation}
\label{eqn:zomega}
E_{-n}\equiv \pi(s_{1}^{n})\Omega,\quad
E_{0}\equiv \Omega,\quad
E_{n}\equiv\pi(s(z^{(1)})\cdots s(z^{(n)}))^{*}\Omega\quad(n\geq 1).
\end{equation}
Then we see that $\pi(s(z^{(n)}))E_{n}=E_{n-1}$ for each $n\in {\bf Z}$
and $\{E_{n}\}_{n\in{\bf Z}}$ is an orthonormal family in ${\cal H}$.
Define
\[E_{n,j}\equiv \pi(s(z^{(n+1,j)}))E_{n+1}\quad((n,j)\in{\bf Z}\times\nset{})\]
where $z^{(n,j)}\equiv(\,g_{j1}(n),\ldots,g_{jN}(n)\,)\in \ct$ 
for $n\in {\bf Z}$ and $j=1,\ldots,N$.
Then $\{E_{n,j}:(n,j)\in{\bf Z}\times\nset{}\}$ is an orthonormal family
in ${\cal H}$.
Especially, $E_{n,1}=E_{n}$ for each $n\in {\bf Z}$.
Define
\[E_{J,n,j}\equiv\pi(s_{J})E_{n+m,j}\quad
(J\in \nset{m},\,n\in{\bf Z},\,m\in {\bf N},\,j=2,\ldots,N).\]
By the estimation of values of inner products,
we see that 
$\{E_{n}\}_{n\in{\bf Z}}\cup\{E_{J,n,j}:J\in\nset{*},\,n\in {\bf Z},\,
j=2,\ldots,N\}$
is an orthonormal family in ${\cal H}$.
For $z\in \sci$, define
\[\mbox{$
\Lambda(z)\equiv \coprod_{(n,m)\in{\bf Z}\times {\bf N}}
\Lambda_{m}^{(n)}(z)$},\]
$\Lambda_{1}^{(n)}(z)\equiv \{(n,0)\}$,
$\Lambda_{2}^{(n)}(z)\equiv\{(n,z^{(n+1,j)})\}_{j=2}^{N}$,
$\Lambda_{m}^{(n)}(z)\equiv\{(n,\vep_{J}\otimes z^{(n+m-1,j)}):j=2,\ldots,N,\,
J\in \nset{m-2}\}$ for $m\geq 3$.
By a proof of similarity to that of Proposition \ref{prop:main},
the following holds.
%
% Proposition 2.3
%
\begin{prop}
\label{prop:basis}
For $z\in \sci$,
let $({\cal H},\pi)$ be $GP(z)$ with the GP vector $\Omega$
and let $\{E_{n}\}_{n\in {\bf Z}}$ be as in (\ref{eqn:zomega}).
Define
\[E_{(n,y)}\equiv \pi(s(y))E_{n+m-1}\quad
((n,y)\in \Lambda^{(n)}_{m}(z),\,(n,m)\in{\bf Z}\times {\bf N}).\]
Then $\{E_{x}:x\in \Lambda(z)\}$ is a complete orthonormal basis of ${\cal H}$.
\end{prop}

\noindent
{\it Proof of Theorem \ref{Thm:mainfirst}.}
\noindent
By Proposition \ref{prop:existence},
the existence is shown.
Assume that both $({\cal H},\pi)$ and $({\cal H}^{'},\pi^{'})$ are $GP(z)$.
According to the natural correspondence among bases of ${\cal H}$
and ${\cal H}^{'}$ in Proposition \ref{prop:main} and \ref{prop:basis},
we have a unitary $U$ from ${\cal H}$ to ${\cal H}^{'}$
such that $U\pi(\cdot)U^{*}=\pi^{'}(\cdot)$.
Hence for $z\in \sc{N}^{\#}$, 
any two GP representations of $\con$ by $z$ are equivalent.
Therefore the uniqueness is proved.
\qedh

%
% Proposition 2.4
%
\begin{prop}
\label{prop:nondeg}
Assume that $z\in \sct{N}{*}$ and 
$({\cal H},\pi)$ is a representation of $\con$
with a cyclic vector $\Omega\in{\cal H}$ such that $\pi(s(z))\Omega=\Omega$,
If $z$ is nonperiodic, then $({\cal H},\pi)$ is $GP(z)$.
\end{prop}
%
% Proof
%
\pr
Assume that $z=z^{(1)}\otimes \cdots \otimes z^{(k)}$ is nonperiodic.
Define $y_{1}\equiv z$, 
$y_{n}\equiv z^{(n)}\otimes \cdots\otimes z^{(k)}\otimes
z^{(1)}\otimes \cdots\otimes z^{(n-1)}$ for $n=2,\ldots,k$
and $E_{n}\equiv \pi(s(z^{(n)})\cdots s(z^{(k)}))\Omega$ for $n=1,\ldots,k$.
Then we have
$<E_{n}|E_{m}>=<\pi(s(y_{n}))E_{n}|\pi(s(y_{m}))E_{m}>=<y_{n}|y_{m}> <E_{n}|E_{m}>$.
By the nonperiodicity of $z$, $|<y_{n}|y_{m}>|<1$ when $n\ne m$.
Hence $|<E_{n}|E_{m}>|=|<y_{n}|y_{m}>|\cdot | <E_{n}|E_{m}>|<| <E_{n}|E_{m}>|$.
This implies that $<E_{n}|E_{m}>=\delta_{n,m}$ for $n,m=1,\ldots,k$.
Hence the statement holds.
\qedh

%%%%%%%%%%%%%%%%%%%%%%%%%%%%%%%%%%%%%%%%%%%%
%
% Section 3
%
\sftt{States and induced product representations}
\label{section:third}
In the first, 
we show the relation between GP representations and states of $\con$.
Next, we review results in \cite{BC}.
At last, we show a relation between
GP representations with a chain and induced product representations.

For $J\in\nset{*}$, the {\it length} $|J|$ of $J$ is defined by
$|J|=k$ when $J\in\nset{k}$.
For $z\in \sc{N}^{\#}$ and $J=(j_{1},\ldots,j_{m})\in\{1,\ldots,N\}^{m}$,
define $z(J)\equiv z^{(1)}_{j_{1}}\cdots z^{(m)}_{j_{m}}$
and $z(0)\equiv 1$ where we define $z^{(kn+l)}_{i}\equiv z^{(l)}_{i}$
for each $n\in {\bf N}$, $l=1,\ldots,k$ and $i\edot$
when $z=z^{(1)}\otimes \cdots\otimes z^{(k)}$.

%
% Proposition 3.1
%
\begin{prop}
\label{prop:statefirst}
For $z\in \sc{N}^{\#}$, define a state $\omega_{z}$ of $\con$ as follows.
When $z=z^{(1)}\otimes\cdots\otimes z^{(k)}\in \ctk$,
\[\omega_{z}(s_{J}s_{K}^{*})\equiv\overline{z(J)}\cdot z(K)\quad(|J|\equiv |K| 
\mbox{ mod }k),\quad \omega_{z}(s_{J}s_{K}^{*})\equiv 0\quad(\mbox{otherwise})\]
for each $J,K\in \nset{*}$
with the convention $s_{J}s_{K}^{*}\equiv s_{K}^{*}$ when $J=\emptyset$.
When $z=(z^{(n)})_{n\in {\bf N}}\in\sci$,
%
% Equation 3.1
%
\begin{equation}
\label{eqn:chainstate}
\omega_{z}(s_{J}s_{K}^{*})\equiv\delta_{|J|,|K|}\cdot \overline{z(J)}\cdot
z(K)\quad(J,K\in \nset{*}).
\end{equation}
Then the GNS representation by $\omega_{z}$ is $GP(z)$.
\end{prop}
%
% Proof
%
\pr
Let $({\cal H},\pi)$ be $GP(z)$ with the GP vector $\Omega$.
Define $\omega^{'}(x)\equiv<\Omega|\pi(x)\Omega>$ for $x\in \con$.
Then we see that $\omega^{'}=\omega_{z}$.
By the uniqueness of GNS representation and the cyclicity of GP representation, 
the statement holds.
\qedh

Let ${\cal H}$ be a Hilbert space in $\con$ defined by 
${\cal H}\equiv {\rm Lin}\langle\{s_{1},\ldots,s_{N}\}\rangle$
with the inner product  $<s_{i}|s_{j}>\equiv \delta_{ij}$ for $i,j\edot$.
For $z=(z^{(i)})\in \sc{N}^{\infty}$,
define a map $s_{r}:{\cal H}^{\otimes r}\to {\cal H}^{\otimes (r+1)}$
by $s_{r}(\xi)\equiv \xi\otimes z^{(r+1)}$.
Define the inductive limit ${\cal H}_{(0)}^{z}$ 
of the inductive system $\{({\cal H}^{\otimes r},s_{r})\}_{r\geq 1}$
of Hilbert spaces
and a vector $\Omega_{(0)}^{z}$ in ${\cal H}_{0}^{z}$ by
$\Omega_{(0)}^{z}\equiv\lim_{r}z^{(1)}\otimes \cdots\otimes z^{(r)}$.
Define subalgebras ${\cal B}_{n}\equiv C^{*}\langle\{s_{J}s_{K}^{*}:
J,K\in \nset{n}\}\rangle$ and
$\con^{0}\equiv \overline{\bigcup_{n=1}^{\infty}{\cal B}_{n}}$ of $\con$,
and denote
$v_{L}\equiv s_{l_{1}}\otimes\cdots\otimes s_{l_{n}}\in {\cal H}^{\otimes n}$
for $L=(l_{1},\ldots,l_{n})\in\nset{n}$.
Remark that $\con^{0}$ equals to the fixed-point subalgebra $\con^{U(1)}$
of $\con$ by the $U(1)$-gauge action.
Define a representation $\pi_{n}$ of ${\cal B}_{n}$ on ${\cal H}^{\otimes n}$ by
\[\pi_{n}(s_{J}s_{K}^{*})v_{L}\equiv \delta_{KL}\cdot v_{J}\quad
(J,K,L\in\nset{n}).\]
Then we have the inductive limit representation $a_{(0)}^{z}$ of $\con^{U(1)}$ 
on ${\cal H}_{(0)}^{z}$ by $\{({\cal H}^{\otimes n},\pi_{n})\}_{n\geq 1}$.
Define a state $F_{(0)}^{z}\equiv
<\Omega_{(0)}^{z}|a_{(0)}^{z}(\cdot)\Omega_{(0)}^{z}>$.
%
% Definition 3.2
%
\begin{defi}
(Definition 2.9, \cite{BC})
\label{defi:ipr}
$({\cal H}^{z},a^{z})$ is the induced product representation of $\con$ by $z$ 
if $({\cal H}^{z},a^{z},\Omega^{z})$ is
the GNS representation of $\con$ of $F^{z}$
where $F^{z}$ is a state of $\con$ defined by $F_{(0)}^{z}\circ P$
and $P$ is the canonical conditional expectation from $\con$ onto $\con^{U(1)}$.
\end{defi}
In \cite{BC}, the induced product representation is parametrized by
the twosided infinite sequence in ${\cal H}$.
We reformulate that by the onesided sequence here.
%
% Theorem 3.3
%
\begin{Thm}
\label{Thm:iprone}
For $z\in\sci$, the following holds.
\begin{enumerate}
%(i)
\item
$a^{z}\sim a^{y}$ if and only if there is $k\geq 0$ such that 
$\sum_{n=1}^{\infty}(1-|<z^{(n)}|y^{(n+k)}>|)<\infty$ or
$\sum_{n=1}^{\infty}(1-|<z^{(n+k)}|y^{(n)}>|)<\infty$.
%(ii)
\item
$a^{z}$ is irreducible if and only if 
$\sum_{n=1}^{\infty}(1-|<z^{(n)}|z^{(n+k)}>|)=\infty$ for any $k\geq 1$.
%(iii)
\item
There is a family $\{a_{n}^{z}\}_{n\in{\bf Z}}$ of (nonzero) subrepresentations
of $a^{z}|_{\con^{U(1)}}$ such that
$a^{z}|_{\con^{U(1)}}=\bigoplus_{n\in {\bf Z}}a^{z}_{n}$.
%(iv)
\item
$a^{z}$ is $GP(z)$.
\end{enumerate}
\end{Thm}
%
% Proof
%
\pr
(i), (ii) and (iii) are shown in Theorem 3.16, 
Corollary 3.17 and Proposition 2.12 in \cite{BC}, respectively.
For states $F^{z}$ in Definition \ref{defi:ipr}
and $\omega_{z}$ in (\ref{eqn:chainstate}),
we can verify that $F^{z}=\omega_{z}$.
This implies (iv).
\qedh

%%%%%%%%%%%%%%%%%%%%%%%%%%%%%%%%%%%%%%%%%
%
% section 4
%
\sftt{Irreducibility, equivalence and decomposition}
\label{section:fourth}
%%%%%%%%%%%%%%%%%%%%%%%%%%%%%%%%%%%%%%%%%%%%%%%%%%
%%%%%%%%%%%%%%%%%%%%%%%%%%%%%%%%%%%%%%%%%%%%%%
%
% subsection 4.1
%
\ssft{Irreducibility}
\label{subsection:fourthone}
Let $\sct{N}{*}_{P}$, $\sci_{EP}$, $\sci_{AP}$ be sets of 
all periodic elements in $\sct{N}{*}$, all \evp\ elements in $\sci$,
all \asp\ elements in $\sci$, respectively, and let 
$\sct{N}{*}_{NP}$, $\sci_{NEP}$, $\sci_{NAP}$ be their complements, respectively.
Then we see that 
$\sct{N}{*}_{P}=\{v^{\otimes k}:v\in\sct{N}{*}_{NP},\, k\geq 2\}$.
%
% Lemma 4.1
%
\begin{lem}
\label{lem:van}
For $z\in \sct{N}{*}_{NP}$, let $({\cal H},\pi)$ be $GP(z)$ 
with the GP vector $\Omega$ and let
$\{E_{x}:x\in \Lambda(z)\}$ be the basis of ${\cal H}$ in (\ref{eqn:basis}).
Then the following holds.
\begin{enumerate}
%(i)
\item
If $z\in \ctk$ and $b\in {\bf N}$ is not a multiple of $k$,
then $\{\pi(s(z)^{*})\}^{m}E_{x}$ goes to $0$
when $m\to \infty$ for any $x\in \Lambda(z)\cap \ct^{\otimes b}$.
%(ii)
\item
If $v\in {\cal H}$ satisfies that $<v|\Omega>=0$, then
$\lim_{m\to \infty}\{\pi(s(z)^{*})\}^{m}v=0$.
\end{enumerate}
\end{lem}
%
% Proof
%
\pr
(i)
For sufficient large $m\in {\bf N}$,
there are $c_{m}\in {\bf C}$ and $j\in\{1,\ldots,k\}$
such that $|c_{m}|\leq 1$ and $\{\pi(s(z)^{*})\}^{m}E_{x}=c_{m}\cdot E_{j}$.
By the assumption on $b$, we see that $j\ne 1$.
The nonperiodicity of $z$ implies that $c_{m}\to 0$ when $m\to \infty$.

\noindent
(ii)
Because $v$ is written as $\sum_{x\in\Lambda(z)}a_{x}E_{x}$,
it is sufficient to consider the case $v=\pi(s(x))\Omega$
for $x\in\Lambda(z)$.
Assume that $z\in \ctk$ and
$v=\pi(s(x))\Omega$ for $x\in \ct^{\otimes b}\cap \Lambda(z)$.
When $b$ is not a multiple of $k$,  the statement holds by (i).
Assume that $b=kl$ for $l\geq 1$.
By assumption, $x\ne z$.
By checking every element $x\in (\Lambda(z)\setminus \{z\})\cap \sct{N}{b}$,
we can verify that $\pi(s(z)^{*})^{m}v= 0$ when $m\geq l$.
In consequence, the assertion holds.
\qedh

%
% Corollary 4.2
%
\begin{cor}
\label{cor:unii}
Let $z\in S({\bf C}^{N})_{NP}^{\otimes k}$ and 
let $({\cal H},\pi)$ be a representation of $\con$.
If $\Omega,\Omega^{'}\in{\cal H}$ are cyclic vectors which 
satisfy that $\pi(s(z))\Omega=\Omega$ and $\pi(s(z))\Omega^{'}=\Omega^{'}$,
then there is $c\in{\bf C}$ such that $\Omega=c\Omega^{'}$.
\end{cor}
%
% Proof
%
\pr
There are $y\in {\cal H}$ and $c\in {\bf C}$ such that 
$\Omega^{'}=c\Omega+y$ and $<\Omega|y>=0$.
By assumption,
we see that $\{\pi(s(z)^{*})\}^{n}\Omega^{'}=\Omega^{'}$ for each $n\in {\bf N}$.
By Lemma \ref{lem:van} (ii), $\lim_{n\to\infty}\{\pi(s(z)^{*})\}^{n}y=0$.
Hence $\Omega^{'}=\lim_{n\to \infty}\{\pi(s(z)^{*})\}^{n}\Omega^{'}$ $=c\Omega$.
\qedh

%
% Proposition 4.3
%
\begin{prop}
\label{prop:irr}
\begin{enumerate}
%(i)
\item
If $z\in \sct{N}{*}_{NP}$, then $GP(z)$ is irreducible.
%(ii)
\item
For $z\in \sci$, $GP(z)$ is irreducible if and only if $z\in \sci_{NAP}$.
\end{enumerate}
\end{prop}
%
% Proof
%
\pr
(i)
Let $({\cal H},\pi)$ be $GP(z)$ with the GP vector $\Omega$.
For a nonzero vector $v_{0}\in {\cal H}$,
it is sufficient to show that  $\Omega\in \overline{\pi(\con)v_{0}}$.
Because ${\cal H}$ is cyclic,
there is $x\in\con$ such that $<\pi(x)\Omega|v_{0}>\ne 0$.
Therefore we can assume that $<\Omega|v_{0}>=1$ by replacing
$v$ and $\pi(x^{*})v$ and normalizing it.
Hence we can denote
$v_{0}=\Omega+y$ for $y\in{\cal H}$ such that $<\Omega|y>=0$.
Assume that $z\in \sct{N}{k}_{NP}$.
By Lemma \ref{lem:van},
we see that $\lim_{n\to \infty}(\pi(s(z)^{*}))^{n}v_{0}=\Omega$.
Hence $\Omega\in \overline{\pi(\con) v_{0}}$.
Therefore $({\cal H},\pi)$ is irreducible.

\noindent
(ii)
This holds by Theorem \ref{Thm:iprone} (ii) and (iv).
\qedh

\noindent
The inverse of Proposition \ref{prop:irr} (i)
is shown in $\S$ \ref{subsection:fourththree}.

%%%%%%%%%%%%%%%%%%%%%%%%%%%%%%%%%%%%%
%
% Subsection 4.2
%
\ssft{Equivalence}
\label{subsection:fourthtwo}
%
%  Lemma 4.4
%
\begin{lem}
\label{lem:pone}
Let $z,y\in \sct{N}{*}$ and let 
$({\cal H},\pi)$ be a representation of $\con$.
If $z\not\sim y$ and there are $\Omega,\Omega^{'}\in{\cal H}$
which satisfy $\pi(s(z))\Omega=\Omega$ and $\pi(s(y))\Omega^{'}=\Omega^{'}$,
then $<\Omega|\Omega^{'}>=0$.
\end{lem}
%
% Proof
%
\pr
Assume that $z\in \ctk$ and $y\in \ct^{\otimes l}$.
If $k\ne l$, then
$<\Omega|\Omega^{'}>=<\pi(s(z)^{l})\Omega|\pi(s(y)^{k})\Omega^{'}>
=<z^{\otimes l}|y^{\otimes k}><\Omega|\Omega^{'}>$.
Because $z\not\sim y$, $|<z^{\otimes l}|y^{\otimes k}>|<1$.
Hence $<\Omega|\Omega^{'}>=0$.
Assume that $k=l$.
If $|<z|y>|<1$, then we see that $<\Omega|\Omega^{'}>=0$.
If $|<z|y>|=1$, then there is $c\in U(1)$ such that $y=cz$.
Since $z\not\sim y$, $c\ne 1$.
By $\pi(s(y))\Omega=c\pi(s(z))\Omega=c\Omega$, $<\Omega|\Omega^{'}>=0$.
\qedh

%
% Proposition 4.5
% 
\begin{prop}
\label{prop:eed}
If $z,y\in \sct{N}{*}$ or $z,y\in \sci$,
then $GP(z)\sim GP(y)$ if and only if $z\sim y$.
\end{prop}
%
% Proof
%
\pr
Let $z,y\in \sct{N}{*}$.
Assume that $z\sim y$ and $z\in \ctk$.
Let $({\cal H},\pi)$ be $GP(z)$ with the GP vector $\Omega$.
By assumption, there is $n\in\{1,\ldots,k\}$ such that 
$y=z^{(n)}\otimes \cdots\otimes z^{(k)}\otimes z^{(1)}\otimes
\cdots\otimes z^{(n-1)}$.
Then $E_{n}\equiv \pi(s(z^{(n)}\otimes \cdots\otimes z^{(k)}))\Omega$
satisfies $\pi(s(y))E_{n}=E_{n}$ and $E_{n}$ is also a  cyclic vector.
Hence $({\cal H},\pi)$ is $GP(y)$.
By the second statement in Theorem \ref{Thm:mainfirst}, $GP(z)\sim GP(y)$.

Assume by contradiction that $z\not \sim y$ and $GP(z)\sim GP(y)$.
By assumption, there is a representation $({\cal H},\pi)$ of $\con$
such that there are unit cyclic vectors $\Omega,\Omega^{'}\in {\cal H}$
which satisfy $\pi(s(z))\Omega=\Omega$ and $\pi(s(y))\Omega^{'}=\Omega^{'}$.
By the basis of $GP(z)$ in (\ref{eqn:basis}),
there are $a_{1},\ldots,a_{k}\in {\bf C}$ and $v\in {\cal H}$
such that $\Omega^{'}=a_{1}E_{1}+\cdots +a_{k}E_{k}+v$
where $E_{1}=\Omega$ and $<E_{n}|v>=0$ for each $n=1,\ldots,k$.
By Lemma \ref{lem:pone} and $z\not\sim y$,
we see that $<\Omega^{'}|E_{n}>=0$ for each $n=1,\ldots,k$.
Hence $\Omega\in V^{\perp}$ for 
$V\equiv{\rm Lin}\langle\{E_{n}\}_{n=1}^{k}\rangle$.
By definition of $\{E_{x}\}_{x\in \Lambda(z)}$,
there is the smallest $m_{0}\in {\bf N}$ such that
$\Omega^{'}\in (V_{m_{0}})^{\perp}$ where
$V_{m_{0}}\equiv{\rm Lin}\langle\{E_{x}:x\in \Lambda^{(m)}_{n}(z),\,
m\leq m_{0},\,n=1,\ldots,k\}\rangle$.
On the other hand,
we see that $\Omega^{'}=\pi(s(y))\Omega^{'}\in(V_{m_{0}+1})^{\perp}$
by definition of $\{E_{x}\}_{x\in \Lambda(z)}$.
This contradicts the choice of $m_{0}$.
Hence $GP(z)\not\sim GP(y)$.

If $z,y\in \sci$,
the statement holds by Theorem \ref{Thm:iprone} (i) and (iv).
\qedh

We consider a relation between chain and cycle here.
For a sequence $\zeta\equiv (z^{(l)})_{l=1}^{k}$ in $\sc{N}$,
define $\zeta^{\infty}\equiv (z^{(n)})_{n\in {\bf N}}\in \sci$
by $z^{(nk+l)}\equiv z^{(i)}$ for each $n\geq 0$ and $l=1,\ldots,k$.
For $\xi=(y^{(l)})_{l=1}^{k}$,
if $z\equiv z^{(1)}\otimes\cdots \otimes z^{(k)}
=y^{(1)}\otimes\cdots \otimes y^{(k)}$,
then $\zeta^{\infty}\sim \xi^{\infty}$.
We denote $\zeta^{\infty}$ by $z^{\infty}$ when there is no ambiguity.
For $z\in \sct{N}{*}$,
the symbol $z^{\infty}$ makes sense as an equivalence class in $\sci$.
This notation is convenient to describe decomposition of representations
in $\S$ \ref{subsection:fourththree}.
By definition, the following holds immediately.
%
% Lemma 4.6
%
\begin{lem}
\label{lem:regular}
\begin{enumerate}
%(i)
\item
If $z\in \sci_{EP}$, then there is $y \in \sct{N}{*}_{NP}$
such that $z\sim y^{\infty}$.
%(ii)
\item
If $z,y\in \sct{N}{*}_{NP}$ satisfy that $z^{\infty}\sim y^{\infty}$,
then there is $c\in U(1)$ such that $z\sim c y$.
\end{enumerate}
\end{lem}

In order to prove Theorem \ref{Thm:mainthird},
we show a branching law of the restriction of a representation
of $\con$ on $\con^{U(1)}$.
%
% Lemma 4.7
%
\begin{lem}
\label{lem:branching}
Let $({\cal H},\pi)$ be $GP(z)$ of $\con$ for 
$z=z^{(1)}\otimes\cdots\otimes z^{(k)}\in\sct{N}{k}_{NP}$ 
with the GP vector $\Omega$.
Then the following irreducible decomposition holds.
%
% Equation 4.1
%
\begin{equation}
\label{eqn:cycle}
({\cal H},\pi|_{\con^{U(1)}})=(V_{1},\pi|_{\con^{U(1)}})
\oplus\cdots\oplus(V_{k},\pi|_{\con^{U(1)}})
\end{equation}
where $V_{i}$ is the completion of $V_{i,0}\equiv\pi(\con^{U(1)})e_{i}$ 
and $e_{i}\equiv \pi(s(z^{(i)})\cdots s(z^{(k)}))\Omega$ for $i=1,\ldots,k$.
\end{lem}
%
% Proof
%
\pr
Here we denote $\pi(s_{i})$ by $s_{i}$ simply.
Define $E_{JK}\equiv s_{J}s_{K}^{*}$ for $J,K\in\nset{l}$.
Then we see that
$<E_{JK}e_{i}|E_{J^{'}K^{'}}e_{j}>=0$ when $i\ne j$
for each $J,K\in\nset{l}$ and $J^{'},K^{'}\in\nset{l^{'}}$.
Therefore $V_{i}$ and $V_{j}$ are orthogonal when $i\ne j$.
On the other hand,
if $|J|=kn+i$ for $n\geq 0$ and $i=0,1,\ldots,k-1$,
then
let $v\equiv z^{(k-i+1)}\otimes\cdots \otimes z^{(k)}\otimes z^{\otimes n}$.
Then $s_{J}\Omega=s_{J}s(v)^{*}e_{k-i+1}$.
Because $s_{J}s(v)^{*}\in \con^{U(1)}$, $s_{J}\Omega\in V_{i}$.
Because ${\rm Lin}\langle\{s_{J}\Omega:J\in\nset{*}\}\rangle$
is dense in ${\cal H}$,
(\ref{eqn:cycle}) holds as a $\con^{U(1)}$-module.

Next we show the irreducibility of $V_{j}$.
If $x\in V_{j}$, then there is $J^{'}\in \nset{kl}$ such that
$c\equiv <e_{j}|s_{J^{'}}^{*}x>\ne 0$.
We replace $x$ by $c^{-1}\cdot s_{J^{'}}^{*}x$.
Then we have a decomposition $x=e_{j}+y$ where $<y|e_{j}>=0$.
Define $X\equiv s(z^{(j)})\cdots s(z^{(k)})s(z^{(1)})\cdots s(z^{(j-1)})$
and $T_{n}\equiv X^{n}(X^{*})^{n}$ for $n\in {\bf N}$. 
Then $T_{n}\in \con^{U(1)}$ and
$T_{n}x\to e_{j}$ when $n\to\infty$ because $z$ is nonperiodic.
Hence $e_{j}\in \overline{\con^{U(1)}x}=V_{j}$.
Therefore $V_{j}$ is an irreducible $\con^{U(1)}$-module for each $j$.
Hence the statement holds.
\qedh

\noindent
{\it Proof of Theorem \ref{Thm:mainthird}.}
By Proposition \ref{prop:eed},
it is sufficient to show that 
for any $y\in \sct{N}{*}$ and $z\in \sci$, $GP(y)\not\sim GP(z)$.
By Lemma \ref{lem:branching} and Theorem \ref{Thm:iprone} (iii),
branching numbers of components of $GP(y)|_{\con^{U(1)}}$ and $GP(z)|_{\con^{U(1)}}$
are always finite and infinite, respectively.
Hence $GP(y)|_{\con^{U(1)}}$ and $GP(z)|_{\con^{U(1)}}$ are never equivalent.
From this, $GP(y)$ and $GP(z)$ are also never equivalent.
\qedh

%%%%%%%%%%%%%%%%%%%%%%%%%%%%%%%%%%%%%%%%%%%%%%%%%%%%%%%
%
%  Section 4.3
%
\ssft{Decomposition}
\label{subsection:fourththree}
%%%%%%%%%%%%%%%%%%%%%%%%%%%%%%%%%%%%%%%%%%%%%%%%%%%%%%%%%%%%%%
%
% Lemma 4.8
%
\begin{lem}
\label{lem:dis}
Let $({\cal H},\pi)$ be a representation of $\con$.
\begin{enumerate}
%(i)
\item
Assume that there are $\Omega_{1},\ldots,\Omega_{M}\in {\cal H}$
and $z_{1},\ldots,z_{M}\in \sct{N}{*}_{NP}$
such that $\pi(s(z_{i}))\Omega_{i}=\Omega_{i}$ for each $i=1,\ldots,M$.
If $z_{i}\not\sim z_{j}$ when $i\ne j$,
then there is a subrepresentation $(V,\pi|_{V})$ of $({\cal H},\pi)$ 
such that $(V,\pi|_{V})\sim GP(z_{1})\oplus\cdots\oplus GP(z_{M})$.
%(ii)
\item
If there are $z_{1},\ldots,z_{M}$,
$y_{1},\ldots,y_{M^{'}}\in \sct{N}{*}_{NP}$ such that
$({\cal H},\pi)\sim\bigoplus_{j=1}^{M}GP(z_{j})$
and $({\cal H},\pi)\sim\bigoplus_{l=1}^{M^{'}}GP(y_{l})$,
then $M=M^{'}$ and there is $\sigma\in {\goth S}_{M}$
such that $z_{\sigma(i)}\sim y_{i}$ for each $i=1,\ldots,M$.
\end{enumerate}
\end{lem}
%
% Proof
%
\pr
(i)
If $M=1$, then it holds by definition.
Assume that $M\geq 2$.
Let $V_{i}\equiv \overline{\pi(\con)\Omega_{i}}$.
Then $(V_{i},\pi|_{V_{i}})$ is $GP(z_{i})$.
Since $z_{i}$ is nonperiodic, $GP(z_{i})$ is irreducible by
Proposition \ref{prop:irr} (i).
By assumption and Proposition \ref{prop:eed},
$\pi|_{V_{i}}\not\sim\pi|_{V_{j}}$ when $i\ne j$.
Hence a subspace $V\equiv\bigoplus_{i=1}^{M}V_{i}$
satisfies the condition in the statement

\noindent
(ii)
By Proposition \ref{prop:irr} (i),
both $GP(z_{i})$ and $GP(y_{j})$ are irreducible for each $i=1,\ldots,M$
and $j=1,\ldots,M^{'}$.
This implies that $M=M^{'}$ and there is $\sigma\in {\goth S}_{M}$
such that $GP(y_{i})\sim GP(z_{\sigma(i)})$
for each $i=1,\ldots,M$.
By Proposition \ref{prop:eed}, $z_{\sigma(i)}\sim y_{i}$
for each $i=1,\ldots,M$.
\qedh

%
% Lemma 4.9
% 
\begin{lem}
\label{lem:ortho}
For $v\in \sct{N}{*}$ and $p\geq 2$,
let $({\cal H},\pi)$ be $GP(v^{\otimes p})$ with the GP vector $\Omega$.
Then $W\equiv {\rm Lin}\langle\{\pi(s(v^{\otimes j}))\Omega\}_{j=1}^{p}\rangle$
is spanned by eigenvectors $\Omega_{1},\ldots,\Omega_{p}$ of $\pi(s(v))$
with eigenvalues $1,e^{2\pi\sqrt{-1}/p},\ldots,e^{2\pi \sqrt{-1}(p-1)/p}$.
\end{lem}
%
% Proof
%	
\pr
By definition, ${\rm dim}W=p$ and $A\equiv \pi(s(v))|_{W}$
is an action of ${\bf Z}_{p}$ on $W$.
Therefore $W$ is decomposed into subspaces $W_{1},\ldots,W_{p}$ as
the irreducible decomposition.
Choose $h_{i}\in W_{i}$ such that $h_{i}\ne 0$.
Then there is $\xi_{i}\in {\bf C}$ such that $Ah_{i}=\xi_{i}h_{i}$
for each $i=1,\ldots,p$.
A representation $W$ of ${\bf Z}_{p}$ is cyclic if and only if
any component in the irreducible decomposition of $W$ 
has multiplicity $1$.
Hence $\{\xi_{l}\}_{l=1}^{p}=\{e^{2\pi\sqrt{-1}l/p}\}_{l=1}^{p}$.
\qedh
 
%
% Proposition 4.10
% 
\begin{prop}
\label{prop:decofi}
If $({\cal H},\pi)$ is $GP(v^{\otimes p})$
for $p\geq 2$ and $v\in \sct{N}{*}_{NP}$, then
\[({\cal H},\pi)\sim GP(v)\oplus GP(e^{2\pi\sqrt{-1}/p}v)\oplus\cdots\oplus 
GP(e^{2\pi\sqrt{-1}(p-1)/p}v).\]
This decomposition is unique up to unitary equivalence.
\end{prop}
%
% Proof
%
\pr
By Lemma \ref{lem:ortho}, there are normalized eigenvectors 
$\Omega_{1},\ldots,\Omega_{p}$ of $\pi(s(v))$ in ${\cal H}$ 
with eigenvalues $1,e^{2\pi\sqrt{-1}/p},\ldots,e^{2\pi \sqrt{-1}(p-1)/p}$,
respectively.
Define $z_{i}\equiv e^{-2\pi\sqrt{-1}(i-1)/p}v$.
Then we have $\pi(s(z_{i}))\Omega_{i}=\Omega_{i}$
for $i=1,\ldots,p$.
By Lemma \ref{lem:dis} (i), there is a subrepresentation
$(V,\pi|_{V})$ of $({\cal H},\pi)$ such that
%
% Equation 4.2
%
\begin{equation}
\label{eqn:set}
(V,\pi|_{V})\sim GP(v)\oplus GP(e^{2\pi\sqrt{-1}/p} v)\oplus\cdots\oplus
GP(e^{2\pi\sqrt{-1}(p-1)/p} v).
\end{equation}
On the other hand,
the GP vector $\Omega$ of $({\cal H},\pi)$ belongs to
${\rm Lin}\langle\{\Omega_{i}\}_{i=1}^{p}\rangle\subset V$
by definition of $\Omega_{i}$ in Lemma \ref{lem:ortho}.
Because $\Omega$ is a cyclic vector, 
${\cal H}=\overline{\pi(\con)\Omega}\subset V\subset {\cal H}$.
From this and (\ref{eqn:set}), the first statement holds.
The uniqueness holds by Lemma \ref{lem:dis} (ii).
\qedh

\noindent
By Proposition \ref{prop:irr} (i) and Proposition \ref{prop:decofi},
the GP representation of $\con$ with a cycle 
is always decomposed into a finite direct sum of irreducible 
GP representations of $\con$ with a cycle.
\\

\noindent
{\it Proof of Theorem \ref{Thm:mainsecond}.}
(i)
By Proposition \ref{prop:decofi}, 	
if $z\in \sct{N}{*}_{P}$, then $GP(z)$ is not irreducible.
By Proposition \ref{prop:irr} (i), 	
the statement holds.

\noindent
(ii)
This holds by Proposition \ref{prop:irr} (ii).
\qedh

\noindent
{\it Proof of Theorem \ref{Thm:mainfourth}.}
(i) is shown by Proposition \ref{prop:decofi}.
We show that
for $v\in \sct{N}{*}$,
%
% Equation 4.3
% 
\begin{equation}
\label{eqn:chian}
GP(v^{\infty})\sim \int_{U(1)}^{\oplus}GP(cv)\,d\eta(c)
\end{equation}
where $cv$ denotes the scalar multiple of $v$ by $c$.

Let $v=v^{(1)}\otimes \cdots \otimes v^{(k)}\in\sct{N}{k}$,
$v^{(l)}=(v_{1}^{(l)},\ldots,v^{(l)}_{N})\in \sc{N}$ for $l=1,\ldots,k$
and $Y_{\infty}\equiv {\bf Z}\times {\bf N}$.
Choose a sequence $(g(l))_{l=1}^{k}$ in $U(N)$
which satisfies $g(l)_{j1}=v^{(l)}_{j}$ for $j\edot$ and $l=1,\ldots,k$.
Define a representation $\pi$ of $\con$
on ${\cal H}\equiv l_{2}(Y_{\infty})$ by
\[\mbox{$\pi(s_{i})e_{kn+l,m}\equiv 
\sum_{j=1}^{N}(g(l)^{*})_{ji}e_{kn+l-1,N(m-1)+j}\quad
((n,m)\in Y_{\infty},\,l=1,\ldots,k)$}.\]
Because $s(v^{(l)})=\alpha_{g(l)}(s_{1})$,
$\pi(s(v^{(l)}))e_{kn+l,1}=e_{kn+l-1,1}$ and $\pi(s(v))e_{kn,1}=e_{k(n-1),1}$
for each $l=1,\ldots,k$ and $n\in {\bf Z}$.
In the same way of the proof of Proposition \ref{prop:existence},
we see that $({\cal H},\pi)$ is $GP(v^{\infty})$.

Let $Y_{k}\equiv \{1,\ldots,k\}\times {\bf N}$ 
and ${\cal K}\equiv l_{2}(Y_{k})$. 
Define a unitary operator $T$ from ${\cal H}$ to  $L_{2}(U(1),{\cal K})$
by $Te_{kn+l,m}\equiv \zeta_{n}\otimes e_{l+1,m}$ for $l=0,\ldots,k-1$
and $(n,m)\in Y_{\infty}$ where $\zeta_{n}(c)\equiv c^{n}$ for $c\in U(1)$.
Define a representation $\Pi\equiv {\rm Ad}T\circ \pi$ of $\con$
on $L_{2}(U(1),{\cal K})$.
Then
\[\mbox{$\Pi(s_{i})(\zeta_{n}\otimes e_{l,m})=
\sum_{j=1}^{N}(g(l-1)^{*})_{ji}\cdot 
\zeta_{n}\otimes e_{l-1,N(m-1)+j}\quad(l=2,\ldots,k)$},\]
\[\mbox{$\Pi(s_{i})(\zeta_{n}\otimes e_{1,m})=
\sum_{j=1}^{N}(g(k)^{*})_{ji}\cdot\zeta_{n-1}\otimes e_{k,N(m-1)+j}$}\]
for $(n,m)\in Y_{\infty}$ and $i\edot$.
For $c\in U(1)$, define a representation $\pi_{c}$ of $\con$ on ${\cal K}$
by \[\mbox{$\pi_{c}(s_{i})e_{l,m}\equiv \sum_{j=1}^{N}(g(l-1)^{*})_{ji}
\cdot e_{l-1,N(m-1)+j}\quad(l=2,\ldots,k)$},\]
\[\mbox{$\pi_{c}(s_{i})e_{1,m}\equiv \bar{c}\cdot\sum_{j=1}^{N}(g(k)^{*})_{ji}
\cdot e_{k,N(m-1)+j}$}\]
for $i\edot$ and $m\in {\bf N}$.
Then we can verify that $({\cal K},\pi_{c})$ is $GP(cv)$
with the GP vector $e_{1,1}$.
Further $(\Pi(s_{i})\phi)(c)=\pi_{c}(s_{i})\phi(c)$ for each $i\edot$,
$\phi\in L_{2}(U(1),{\cal K})$ and $c\in U(1)$.
Therefore
$GP(v^{\infty})\sim \Pi\sim \int_{U(1)}^{\oplus}
\pi_{c}\,d\eta(c)\sim \int_{U(1)}^{\oplus}GP(cv)\,d\eta(c)$.
Hence (\ref{eqn:chian}) is shown.
From this and Lemma \ref{lem:regular} (i),   
the first statement of (ii) is proved.
Because $GP(v)\circ \gamma_{\bar{c}^{1/k}}=GP(c v)$ for $c\in U(1)$,
(\ref{eqn:chian}) is rewritten as
\[GP(v^{\infty})\sim\int_{U(1)}^{\oplus}GP(v)\circ\gamma_{c^{1/k}}\,d\eta(c)\]
where $\gamma$ is the gauge action on $\con$.
If $z\in \sci_{EP}$, 
then
there is $y\in \scm{k}_{NP}$ such that $GP(z)\sim GP(y^{\infty})$
by Lemma \ref{lem:regular} (i) and Theorem \ref{Thm:mainthird}.
From this,
%
% Equation 4.4
%
\begin{equation}
\label{eqn:deceq}
GP(z)\sim\int_{U(1)}^{\oplus}GP(y)\circ\gamma_{c^{1/k}}\,d\eta(c).
\end{equation}
If there are $l\geq 1$ and $y_{0}\in \scm{l}$
which satisfy (\ref{eqn:deceq}) with respect to $z$,
then 
$GP(y_{0}^{\infty})\sim
\int_{U(1)}^{\oplus}GP(y_{0})\circ\gamma_{c^{1/l}}\,d\eta(c)
\sim GP(z)\sim GP(y^{\infty})$.
By Theorem \ref{Thm:mainthird}, $y_{0}^{\infty}\sim y^{\infty}$.
By Lemma \ref{lem:regular} (ii),
the last statement of (ii) is proved.
\qedh

%%%%%%%%%%%%%%%%%%%%%%%%%%%%%%%%%%%%%%%%%%%%%%%%%
%
% Section 5
%
\sftt{Examples}
\label{section:fifth}
%%%%%%%%%%%%%%%%%%%%%%%%%%%%%%%%%%%%%%%%%%%%%
%%%%%%%%%%%%%%%%%%%%%%%%%%%%%%%%%%%%%%%%%%%%%%%%%%%
%
% subsection 5.1
%
\ssft{Permutative representations}
\label{subsection:fifthone}
We show that GP representations are generalizations of 
permutative representations by \cite{BJ}.
A data $({\cal H},\pi)$ is a {\it permutative representation} of $\con$
if there is a complete orthonormal basis $\{e_{n}\}_{n\in\Lambda}$ 
of ${\cal H}$ and a family $f=\{f_{i}\}_{i=1}^{N}$ of maps
on $\Lambda$ such that $\pi(s_{i})e_{n}=e_{f_{i}(n)}$
for each $n\in\Lambda$ and $i\edot$.
For $J=(j_{l})_{l=1}^{k}$ and $c\in U(1)\equiv \{a\in {\bf C}:|a|=1\}$,
a representation $({\cal H},\pi)$ of $\con$ is $P(J;c)$
if there is a unit cyclic vector $\Omega\in {\cal H}$
and $c\in U(1)$ such that $\pi(s_{J})\Omega=c\Omega$
and $\{\pi(s_{j_{l}}\cdots s_{j_{k}})\Omega\}_{l=1}^{k}$
is an orthonormal family in ${\cal H}$.
We simply denote $P(J;1)$ by $P(J)$.
For $J=(j_{n})_{n\in {\bf N}}\in\nset{\infty}
\equiv \{(i_{n})_{n\in {\bf N}}:i_{n}\in\nset{}\}$,
a representation $({\cal H},\pi)$ of $\con$ is $P(J)$
if there is a unit cyclic vector $\Omega\in {\cal H}$
such that $\{\pi((s_{J_{(n)}})^{*})\Omega:n\in {\bf N}\}$
is an orthonormal family in ${\cal H}$
where $J_{(n)}\equiv (j_{1},\ldots,j_{n})$.
$({\cal H},\pi)$ is a {\it cycle} ({\it chain}) 
if there is $J\in\nset{k}$ ({\it resp.} $J\in\nset{\infty}$)
such that $({\cal H},\pi)$ is $P(J)$.
Then $({\cal H},\pi)$ is a cyclic permutative representation of $\con$
if and only if $({\cal H},\pi)$ is a cycle or a chain.
We see that $P(J;c)=GP(\bar{c}\cdot \vep_{J})$
for $J\in\nset{k}$ and $c\in U(1)$, 
and $P(J)=GP((\vep_{j_{n}})_{n\in {\bf N}})$
for $J=(j_{n})_{n\in{\bf N}}\in \nset{\infty}$.
%
% Example 5.1
%
\begin{ex}{\rm
\begin{enumerate}
% (i)
\item
The {\it standard representation} $(l_{2}({\bf N}),\pi_{S})$ of $\con$
is defined by $\pi_{S}(s_{i})e_{n}\equiv e_{N(n-1)+i}$
for $n\in {\bf N}$ and $i\edot$
where $\{e_{n}\}_{n\in {\bf N}}$ is the canonical basis of $\ltn$ \cite{AK1}.
Then $(l_{2}({\bf N}),\pi_{S})$ is $GP(z)$ with the GP vector $e_{1}$
for $z=(1,0,\ldots,0)\in\ct$.
Since $z$ is nonperiodic,
$(l_{2}({\bf N}),\pi_{S})$ is irreducible.
% (ii)
\item
Define a representation $(\ltn,\pi)$ of $\co{3}$ by
$\pi(s_{1})e_{1}\equiv e_{2},
\pi(s_{1})e_{2}\equiv e_{1},
\pi(s_{2})e_{1}\equiv e_{3},
\pi(s_{3})e_{1}\equiv e_{4}$,
$\pi(s_{2})e_{2}\equiv e_{5},
\pi(s_{3})e_{2}\equiv e_{6},
\pi(s_{i})e_{n}\equiv e_{3(n-1)+i}$ for $i=1,2,3,\,n\geq 3$.
Then $e_{1}$ is a cyclic vector of $(\ltn,\pi)$
and $\pi(s_{1}s_{1})e_{1}=e_{1}$.
Therefore $(\ltn,\pi)$ is $GP((1,0,0)\otimes (1,0,0))$.
We have the irreducible decomposition $\ltn=V_{1}\oplus V_{2}$
where $V_{1}\equiv \overline{\pi(\co{3})(e_{1}+e_{2})}$
and $V_{2}\equiv \overline{\pi(\co{3})(e_{1}-e_{2})}$.
\end{enumerate}
}
\end{ex}
%%%%%%%%%%%%%%%%%%%%%%%%%%%%%%%%%%%%%%%%%%
%
% subsection 5.2
%
\ssft{Non permutative representations}
\label{subsection:fifthtwo}
%
% Example 5.2
%
\begin{ex}
\label{ex:seventhone}
{\rm
Define a representation $(L_{2}[0,1],\pi_{B})$ of $\con$ by 
\[(\pi_{B}(s_{i})\phi)(x)\equiv\chi_{[(i-1)/N,i/N]}(x)\phi(Nx-i+1)\]
for $\phi\in L_{2}[0,1]$, $x\in [0,1]$ and $i\edot$
where $\chi_{Y}$ is the characteristic function of a subset $Y$ of $[0,1]$.
Then $(L_{2}[0,1],\pi_{B})$ is $GP(z)$ for $z\equiv (N^{-1/2},\ldots,N^{-1/2})$.
$(L_{2}[0,1],\pi_{B})$ is also irreducible.
Further $(L_{2}[0,1],\,\pi_{B}\circ \gamma_{\bar{c}})$
is $GP(cz)$ for $c\in U(1)$.
}
\end{ex}

%
% Example 5.3
%
\begin{ex}
{\rm
(i)
For a unitary matrix $g=\mattwo{a}{b}{c}{d}\in U(2)$,
let $(l_{2}({\bf N}),\pi)$ be a representation of $\co{2}$ defined by
\[\pi(s_{1})e_{1}\equiv ae_{2}+ce_{3},\quad\pi(s_{1})e_{2}\equiv e_{1}, \quad
\pi(s_{2})e_{1}\equiv be_{2}+ed_{3},\quad\pi(s_{2})e_{2}\equiv e_{4},\]
\[\pi(s_{1})e_{n}\equiv e_{2n-1},\quad \pi(s_{2})e_{n}\equiv e_{2n}
\quad(n\geq 3).\]
Then $(l_{2}({\bf N}),\pi)$ satisfies that
$\pi(s_{1}(\alpha s_{1}+\beta s_{2}))e_{1}=e_{1}$
for $g^{-1}=\mattwo{\alpha}{\gamma}{\beta}{\delta}\in U(2)$.
Hence this representation is $GP(z)$ for $z\equiv
\vep_{1}\otimes (\alpha\vep_{1}+\beta\vep_{2})$.
If $\alpha\beta\ne 0$,
then this representation does not belong to the class of representation by
\cite{BJ,DaPi2,DaPi3}.	
We see that $\beta\ne 0$ if and only if $z$ is nonperiodic.
This is equivalent that $GP(z)$ is irreducible.
If $\beta=0$, then $|\alpha|=1$ and the following equivalence holds.
\[GP((\alpha^{1/2}\vep_{1})^{\otimes 2})\sim 
GP(\alpha^{1/2}\vep_{1})\oplus GP(-\alpha^{1/2}\vep_{1}).\]
Note that this representation
is determined only by $(\alpha,\beta)\in {\bf C}^{2}$,
$|\alpha|^{2}+|\beta|^{2}=1$ up to unitary equivalence.
Furthermore for any two different elements in $S({\bf C}^{2})$,
the associated representations are inequivalent.
In this sense, $S({\bf C}^{2})$ is a parameter
space of unitary equivalence classes of representations of $\co{2}$
in this case.

\noindent
(ii)
Because any $z\in\ct$ is nonperiodic, a cyclic representation
of $\con$ with the cyclic vector $\Omega$ which satisfies
$\pi(s(z))\Omega=\Omega$ is irreducible by Proposition \ref{prop:irr} (i).
For $z,y\in \sc{N}$, $GP(z)\sim GP(y)$ if and only if $z=y$.

\noindent
(iii)
Let $p\geq 1$ and let $({\cal H},\pi)$ be a representation
of $\con$ with a cyclic vector $\Omega$.
If $({\cal H},\pi)$ satisfies the following condition,
then $({\cal H},\pi)$ is irreducible:
\[\pi((s_{1}+s_{2})(s_{1}+\xi s_{2})(s_{1}+\xi^{2} s_{2})
\cdots(s_{1}+\xi^{p-1} s_{2}))\Omega=2^{p/2}\Omega\]
where $\xi\equiv e^{2\pi\sqrt{-1}/p}$.
}
\end{ex}
%
% Example 5.4
%
\begin{ex}
\label{ex:irrational}
{\rm
For $\theta\in [0,1)$, define 
$z_{\theta}=(z^{(n)}_{\theta})_{n\in {\bf N}}\in S({\bf C}^{2})^{\infty}$ by
\[z^{(n)}_{\theta}\equiv (\cos2\pi n\theta ,\sin2\pi n\theta)
\in S^{1}\equiv \{(x,y)\in {\bf R}^{2}:x^{2}+y^{2}=1\}\quad(n\in {\bf N}).\]
Then $GP(z_{\theta})$ is irreducible if and only if $\theta$ is irrational.
If there is $p_{\theta}\equiv
{\rm min}\{q\in {\bf N}:q\theta \in {\bf N}\cup \{0\}\}$, then
\[GP(z_{\theta})\sim \int_{U(1)}^{\oplus}
GP(c\cdot z_{\theta}^{(1)}\otimes \cdots\otimes z_{\theta}^{(p_{\theta})})
\,d\eta(c).\]
%
% Proof
%
\pr
Because
\[\sum_{n=1}^{M}(1-|<z_{\theta}^{(n)}|z_{\theta}^{(n+p)}>|)
=2M\cdot \sin^{2}\pi p\theta\quad(p,M\geq 1),\]
we see that
there is $p\geq 1$ such that
$\sum_{n=1}^{\infty}(1-|<z_{\theta}^{(n)}|z_{\theta}^{(n+p)}>|)<\infty$
if and only if $\theta\in {\bf Q}$.
Therefore
$z_{\theta}$ is \asp\ if and only if $\theta\in {\bf Q}$.
By Theorem \ref{Thm:mainsecond} (ii), the first statement holds.
If $p_{\theta}$ exists, then $z_{\theta}$ is \evp.
Hence the last statement holds by (\ref{eqn:chian}).
\qedh
}
\end{ex}
%%%%%%%%%%%%%%%%%%%%%%%%%%%%%%%%%%%%%%%%%%%%%%%%%%%%%%%%%%%%%%%%%%%%%
%
% subsection 5.3 
%
\ssft{Eventual periodicity and asymptotic periodicity}
\label{subsection:fifththree}
In order to show a remaining problem, we introduce an example
$z\in S({\bf C}^{2})^{\infty}$ which is \asp\ but not \evp.

%
% Lemma 5.5
% 
\begin{lem}
\label{lem:construction}
Define positive sequences $\{a_{n}\}$ and $\{b_{n}\}$ by
\[c_{n}\equiv \frac{\pi}{4}-Y_{n},\quad d_{n}\equiv \frac{\pi}{4}+Y_{n},\quad 
Y_{n}\equiv \frac{1}{2}\sin^{-1}\frac{1}{\sqrt{2}n}\quad(n\geq 1)\]
and
let $z\equiv (z^{(n)})_{n\geq 1}\in S({\bf C}^{2})^{\infty}$ by
\[z^{(2n-1)}\equiv (\cos c_{n},\sin c_{n}),\quad
z^{(2n)}\equiv (\cos d_{n},\sin d_{n})\quad(n\geq 1).\]
Then the following holds.
\begin{enumerate}
%(i)
\item $z$ is non\evp\ and \asp.
%(ii)
\item
$z$ is equivalent to $(1/\sqrt{2}, 1/\sqrt{2})^{\infty}\in S({\bf C}^{2})^{\infty}$.
\end{enumerate}
\end{lem}
%
% Proof
%
\pr
(i)
Since $z$ is a sequence in ${\bf R}^{2}$ with positive components
and $z^{(n)}\ne z^{(m)}$ when $n\ne m$,
it is non\evp.
We see that $1-|<z^{(2n-1)}|z^{(2n)}>|=1/n^{2}$.
Because $d_{n}-c_{n+1}< d_{n}-c_{n}$,
we can verify that $1-|<z^{(2n)}|z^{(2n+1)}>|<1/n^{2}$.
Hence $\sum_{n=1}^{\infty}(1-|<z^{(n)}|z^{(n+1)}>|)<\pi^{2}/3<\infty$.
Therefore $z$ is \asp.

\noindent
(ii)
Let $v\equiv (1/\sqrt{2}, 1/\sqrt{2})$.
Then we see that
\[1-|<z^{(2n-1)}|v>|=2\sin^{2}(Y_{n}/2)<2\sin^{2}(\sin^{-1}(1/\sqrt{2}n))
=1/n^{2}.\]
As well as this, $1-|<z^{(2n)}|v>|<1/n^{2}$.
In consequence, $\sum_{n=1}^{\infty}(1-|<z^{(n)}|v>|)<\pi^{2}/3$.
Hence the statement holds.
\qedh

\noindent
By Lemma \ref{lem:construction} (i), 
we have the following tabular about chain case:
\[
{\footnotesize
\mbox{
\begin{tabular}{c|c|c}
\hline
\multicolumn{2}{c|}{non\evp}&\evp \\
\hline
indecomposable&gray zone&decomposable\\
\hline
nonasymptotically periodic &\multicolumn{2}{c}{asymptotically periodic}\\
\hline
irreducible&\multicolumn{2}{c}{non irreducible}\\
\hline
\end{tabular}
}
}\]

\noindent
where the decomposability is defined as the style of (\ref{eqn:deceq}).
In the ``gray zone" $=\sci_{NEP}\cap \sci_{AP}$, 
it is not known that whether $GP(z)$ is decomposable or not.
According to Lemma \ref{lem:construction} (ii), we propose the following.
%
% Conjecture 5.6
% 
\begin{conj}
If $z$ is \asp, then there is always $y\in\sci_{EP}$ such that $z\sim y$.
\end{conj}

%%%%%%%%%%%%%%%%%%%%%%%%%%%%%%%%%%%%%%%%%%%%%%%%%%%%%%%%%%%%%%%%%%%%%

\noindent
{\bf Acknowledgement:}
The author would like to thank Prof. Akira Asada 
for a suggestion to study representations of Cuntz algebras.
Florin Ambro and Takeshi Nozawa help us by critical reading of this article.

%%%%%%%%%%%%%%%%%%%%%%%%%%%%%%%%%%%%%%%%%%%%%%%%%%%%%%%%%%%%%%%%%%%%%

\end{document}